\documentclass{article}
\usepackage{amsfonts}
\usepackage{amssymb}
\usepackage{amsmath}
\usepackage{listings}
\usepackage{amssymb}
\usepackage{url}
\usepackage{nicefrac}
\usepackage{dcolumn}
\usepackage{tabularx}
\usepackage{algorithmic}
\usepackage{algorithm}
\newcommand{\ra}[1]{\renewcommand{\arraystretch}{1.2}}
\usepackage{xcolor}
\usepackage[square,sort&compress]{natbib}
\usepackage{bm}
\renewcommand{\maketitle}%
{\thispagestyle{empty}
\noindent
\par\vspace*{3cm}\par
{\centering{\Large\textbf\MATCHtitle\par}
\par\vspace*{6pt}\par
{\textrm\MATCHauthor}\par\vspace*{12pt}\par}%
}
\renewenvironment{abstract}%
{\begin{center}\begin{minipage}[t]{1.0\textwidth}
\noindent\small\textbf{Abstract$\colon$}}%
{\end{minipage}\end{center}\par}

\usepackage{alltt} 
\newcommand{\MATCHtitle}
{How to increase convergence order of Newton's method to $2\times m$?}
\newcommand{\MATCHauthor}
{
}
\begin{document}
\maketitle

\begin{center}
 Sanjay Kumar Khattri \footnote{Corresponding author.\\
E-mail addresses: sanjay.khattri@hsh.no (S. K. Khattri)}\\
Department of Engineering, \\
Stord Haugesund University College, Norway\\
\end{center}

\begin{abstract}
We present a simple yet powerful technique for forming iterative methods of various convergence orders. Methods of various convergence orders (four, six, eight and ten) are formed through a modest modification of the classical Newton method. The technique can be easily implemented in existing software packages as suggested by the presented C$^{++}$ algorithm. Finally some problems are solved through the proposed algorithm.  
\end{abstract}
\section{Introduction}\label{sec:introduction}
The most common and probably the most used method method finding a simple root
$\gamma,$ i.e. $f(\gamma)=0$, of a nonlinear scalar equation
\begin{equation}
f(x) = 0\mathpunct{,}\label{eq:equation}
\end{equation}
is the Newton method. The classical Newton method is given as
\begin{equation}
x_{n+1} = x_{n} - \dfrac{f(x_n)}{f^\prime(x_n)},\qquad n =
0,1,2,3,\ldots. \label{eq:newton}
\end{equation}
It is well documented and well known that the Newton's method converges quadratically \citep[see][and references
therein]{fourth1,fourth2,fourth3,fourth4,fourth5,fourth6,fourth7,fourth8,fourth9,fourth10,abbasbandy,example,fourth11,fourth12,khattri2,khattri_1,thirdorder1,thirdorder2,thirdorder3,thirdorder4, thirdorder5,potra1,ostro,traub,jarrat}. There exists many modifications of the Newton method to improve convergence order \citep{fourth1,fourth2,fourth3,fourth4,fourth5,fourth6,fourth7,fourth8,fourth9,fourth10,abbasbandy,example,fourth11,fourth12,khattri2,khattri_1,thirdorder1,thirdorder2,thirdorder3,thirdorder4, thirdorder5,potra1,ostro,traub,jarrat}. Higher order modifications of the Newton's method free of second or higher derivatives are actively researched. For example, third order convergent methods are presented in \citep{thirdorder2,thirdorder3,thirdorder4, thirdorder5}, fourth order convergent methods are developed in \citep{fourth1,fourth2,fourth3,fourth4,fourth5,fourth6,fourth7,fourth8,fourth9,fourth10}, sixth order methods are developed in  \citep{six1,six2,six3} and eight order methods are presented in \citep[see][and references
therein]{eight1}. As there exists various modifications of the Newtons method. And from practical point (implementing these methods into a software package) one of the drawbacks of these wonderful methods is their independent nature. For example, if one has a software package which solves nonlinear equations by the well known fourth order convergent Jarrat method \citep{jarrat}. One may find it difficult to modify the existing software package to implement sixth order methods \citep{six1,six2,six3} and eight order methods \citep{eight1}.

In this work, we develop a technique. Which improves the order of convergence of the Newton method \eqref{eq:newton} from $2$ to $2\times{m}$. Here, $m=1,2,3,\ldots$. Thus through our scheme, one may develop $4$th order, $6$th order, $8$th order, $\ldots$  convergent iterative methods. One of the beautiful fact of our scheme is that one needs modest modifications in the classical iterative method \eqref{eq:newton} for achieving higher convergence rates. And, which may be very effective when one wants to modify an existing software package for achieving higher convergence order. Let us now develop our scheme.

\section{The technique and convergence order of its various methods}
Before presenting our technique, first we will develop iterative methods of various convergence orders. Let us consider the following $4$th order convergent iterative method
\begin{alignat}{4}
y_n &=& &x_n - \dfrac{f(x_n)}{f^\prime(x_n)}, \label{eq:newtonstep} \\
x_{n+1} &=& &x_n - \dfrac{f(x_n)}{f^\prime(x_n)}\left[1 + \dfrac{f(y_n)}{f(x_n)} \left( 1 + 2\,\dfrac{f(y_n)}{f(x_n)}\right)
\right],\label{eq:fourthorder}
\intertext{error equation for the above method is given as}
e_{n+1} &=& & -\frac{1}{12}\,{\frac {c_{{2}} \left( 12\,c_{{3}}c_{{1}}-60\,{c_{{2}}}^{2}
 \right) }{{c_{{1}}}^{3}}}
e_n^4 + O(e_n^5).
\end{alignat}
Here, $c_k = f^k(\gamma)/{!k}$ and $e_n=x_n-\gamma$. A proof of convergence of the above fourth order method is presented next.
\subsection*{proof}
Using Taylor series of $f(x)$ around the solution $\gamma$, and taking into account $f(\gamma)=0$, we get
\begin{alignat}{4}
f(x_n) &=& &\sum_{k=1}^\infty{c_k\,e_n^k},\label{eq:f(x)}
\intertext{furthermore from the equation \eqref{eq:f(x)} we have}
f^\prime(x_n) &=& &\sum_{k=1}^\infty{k\,c_k\,\,e_n^{k-1}},\label{eq:df(x)}
\intertext{and through a simple calculation we arrive at}
\frac{f(x)}{f^\prime(x)}&=& &e_n-{\frac {c_{{2}}}{c_{{1}}}}{e_n}^{2}-2\,{\frac {c_{{3}}c_{{1}}-{c_{{2}
}}^{2}}{{c_{{1}}}^{2}}}{e_n}^{3}-{\frac {3\,c_{{4}}{c_{{1}}}^{2}-7\,c_{{
2}}c_{{3}}c_{{1}}+4\,{c_{{2}}}^{3}}{{c_{{1}}}^{3}}}{e_n}^{4}+O \left( {e_n
}^{5} \right).\label{eq:fxdfx} 
\intertext{Substituting \eqref{eq:fxdfx} in \eqref{eq:newtonstep} yields}
y_n -\gamma &=& &{\frac {c_{{2}}}{c_{{1}}}}{e_n}^{2}+2\,{\frac {c_{{3}}c_{{1}}+{c_{{2}
}}^{2}}{{c_{{1}}}^{2}}}{e_n}^{3}+{\frac {3\,c_{{4}}{c_{{1}}}^{2}+7\,c_{{
2}}c_{{3}}c_{{1}}+4\,{c_{{2}}}^{3}}{{c_{{1}}}^{3}}}{e_n}^{4}+O \left( {e_n
}^{5} \right).\label{eq:ynminusgamma}
\intertext{Expanding $f(y_n)$ around the solution $\gamma$ and using \eqref{eq:ynminusgamma}, we obtain}
f(y_n) &=& &c_{{2}}{e_n}^{2}-\frac{1}{6}\,{\frac {-12\,c_{{3}}c_{{1}}+12\,{c_{{2}}}^{2}}{c_
{{1}}}}{e_n}^{3}+\frac{1}{24}\,{\frac {72\,c_{{4}}{c_{{1}}}^{2}-168\,c_{{2}}c_{{
3}}c_{{1}}+120\,{c_{{2}}}^{3}}{{c_{{1}}}^{2}}}{e_n}^{4}+O \left( {e_n}^{5}
 \right). \label{eq:f(y)} 
\intertext{From equations \eqref{eq:f(x)} and \eqref{eq:f(y)}, we get}
\frac{f(y_n)}{f(x_n)} &=& &{\frac {c_{{2}}}{c_{{1}}}}e_n+{\frac {2\,c_{{3}}c_{{1}}-3\,{c_{{2}}}^{2
}}{{c_{{1}}}^{2}}}{e_n}^{2}-{\frac {-3\,c_{{4}}{c_{{1}}}^{2}+10\,c_{{2}}
c_{{3}}c_{{1}}-8\,{c_{{2}}}^{3}}{{c_{{1}}}^{3}}}{e_n}^{3}+O \left( {e_n}^{
4} \right). \label{eq:fyfx}
\intertext{Now from equations \eqref{eq:fxdfx}, \eqref{eq:fyfx} and \eqref{eq:fourthorder}, we find that}
e_{n+1} &=& & -\frac{1}{12}\,{\frac {c_{{2}} \left( -60\,{c_{{2}}}^{2}+12\,c_{{3}}c_{{1}}\right) }{{c_{{1}}}^{3}}}e_n^4 + O(e_n^5).
\end{alignat}
This proofs that the method \eqref{eq:fourthorder} converges quartically.

Let us now consider the following three step sixth order convergent iterative method
\begin{eqnarray}
y_n &=& x_n - \dfrac{f(x_n)}{f^\prime(x_n)}, \nonumber \\
z_{n} &=& x_n - \dfrac{f(x_n)}{f^\prime(x_n)}\left[1 + \dfrac{f(y_n)}{f(x_n)}\left(1+  2\,\dfrac{f(y_n)}{f(x_n)}\right)\right].\nonumber \\
x_{n+1} &=& x_n - \dfrac{f(x_n)}{f^\prime(x_n)}\left[1 + \dfrac{f(y_n)}{f(x_n)}\left(1  + 2\,\dfrac{f(y_n)}{f(x_n)}\right)\right. \nonumber \\
& & \left.+ \dfrac{f(z_n)}{f(x_n)} \left(1+2\,\frac{f(y_n)}{f(x_n)}\right)\right],\label{eq:sixthorder}
\end{eqnarray} 
error equations for the above sixth order method is given as
$$e_{n+1} = {\frac {1}{72}}\,{\frac {c_{{2}} \left( -792\,c_{{3}}c_{{1}}{c_{{2}}}^
{2}+2160\,{c_{{2}}}^{4}+72\,{c_{{3}}}^{2}{c_{{1}}}^{2} \right) }{{c_{{
1}}}^{5}}}e_{n}^6+O\left(e_n^7\right).
$$
Convergence order of the above method can be easily established through the Maple software package. We may notice that the method \eqref{eq:sixthorder} requires evaluations of only three functions and one derivative during each iterative step. Let us now further consider the following eight order convergent iterative step
\begin{eqnarray}
y_n &=& x_n - \dfrac{f(x_n)}{f^\prime(x_n)},\nonumber \\
z_{n} &=& x_n - \dfrac{f(x_n)}{f^\prime(x_n)}\left[1 + \dfrac{f(y_n)}{f(x_n)}\left(1+  2\,\dfrac{f(y_n)}{f(x_n)}\right)\right], \nonumber\\
p_{n} &=& x_n - \dfrac{f(x_n)}{f^\prime(x_n)}\left[1 + \dfrac{f(y_n)}{f(x_n)}\left(1 + 2\,\dfrac{f(y_n)}{f(x_n)}\right) + \dfrac{f(z_n)}{f(x_n)} \left(1+2\,\frac{f(y_n)}{f(x_n)}\right)\right],\nonumber\\
x_{n+1} &=& x_n - \dfrac{f(x_n)}{f^\prime(x_n)}\left[1 + \dfrac{f(y_n)}{f(x_n)}\left(1 + 2\,\dfrac{f(y_n)}{f(x_n)}\right)  + \dfrac{f(z_n)}{f(x_n)} \left(1+2\,\frac{f(y_n)}{f(x_n)}\right) \right. \nonumber \\
&& \left.  + \dfrac{f(p_n)}{f(x_n)} \left(1+2\,\frac{f(y_n)}{f(x_n)}\right)\right]\label{eq:eightorder},
\end{eqnarray} 
asymptotic error equation for the eight order method is given as
$$e_{n+1} = {\frac {1}{432}}\,{\frac {c_{{2}} \left( 77760\,{c_{{2}}}^{6}-41472\,c
_{{3}}c_{{1}}{c_{{2}}}^{4}+7344\,{c_{{3}}}^{2}{c_{{1}}}^{2}{c_{{2}}}^{
2}-432\,{c_{{3}}}^{3}{c_{{1}}}^{3} \right) }{{c_{{1}}}^{7}}}e_{n}^8$$
We may notice that the eight order method \eqref{eq:eightorder} requires evaluations of only four functions and one derivative during each iterative step. Based upon the similarity in methods \eqref{eq:fourthorder}, \eqref{eq:sixthorder} and \eqref{eq:eightorder}. Let us consider the following method
\begin{eqnarray}
y_n &=& x_n - \dfrac{f(x_n)}{f^\prime(x_n)},\nonumber \\
z_{n} &=& x_n - \dfrac{f(x_n)}{f^\prime(x_n)}\left[1 + \dfrac{f(y_n)}{f(x_n)} \left(1 + 2\,\dfrac{f(y_n)}{f(x_n)}\right)\right], \nonumber\\
p_{n} &=& x_n - \dfrac{f(x_n)}{f^\prime(x_n)}\left[1 + \dfrac{f(y_n)}{f(x_n)}\left(1 + 2\,\dfrac{f(y_n)}{f(x_n)}\right) + \dfrac{f(z_n)}{f(x_n)} \left(1+2\,\frac{f(y_n)}{f(x_n)}\right)\right],\nonumber\\
q_{n} &=& x_n - \dfrac{f(x_n)}{f^\prime(x_n)}\left[1 + \dfrac{f(y_n)}{f(x_n)}\left(1 + 2\,\dfrac{f(y_n)}{f(x_n)}\right)  + \dfrac{f(z_n)}{f(x_n)} \left(1+2\,\frac{f(y_n)}{f(x_n)}\right)\right.\nonumber\\
&&\left. + \dfrac{f(p_n)}{f(x_n)} \left(1+2\,\frac{f(y_n)}{f(x_n)}\right)\right],\nonumber\\
x_{n+1} &=& x_n - \dfrac{f(x_n)}{f^\prime(x_n)}\left[1 + \dfrac{f(y_n)}{f(x_n)}\left(1 + 2\,\dfrac{f(y_n)}{f(x_n)}\right)  + \dfrac{f(z_n)}{f(x_n)} \left(1+2\,\frac{f(y_n)}{f(x_n)}\right)\right.\nonumber\\
&&\left. + \dfrac{f(p_n)}{f(x_n)} \left(1+2\,\frac{f(y_n)}{f(x_n)}\right)+ \dfrac{f(q_n)}{f(x_n)} \left(1+2\,\frac{f(y_n)}{f(x_n)}\right)\right]\label{eq:tenthorder},
\end{eqnarray}
through the Maple we verified that the above method is $10$th order convergent, and error equation for it is given as
\begin{multline}
e_{n+1} = \frac{c_2}{2592\,c_1^9} \left( 2799360\,{c_{{2}}}^{8}-
1959552\,c_{{3}}c_{{1}}{c_{{2}}}^{6}+513216\,{c_{{3}}}^{2}{c_{{1}}}^{2
}{c_{{2}}}^{4}\right. \nonumber \\
\left. -59616\,{c_{{3}}}^{3}{c_{{1}}}^{3}{c_{{2}}}^{2} +2592\,{c_{{3}}}^{4}{c_{{1}}}^{4} \right) e_n^{10}. \nonumber
\end{multline}
We may notice that the above tenth order method \eqref{eq:tenthorder} requires evaluations of only five functions and one derivative per iterative step.

Based upon the methods \eqref{eq:fourthorder}, \eqref{eq:sixthorder}, \eqref{eq:eightorder} and \eqref{eq:tenthorder}, we conjuncture the existence of the following scheme for generating iterative method of order $2\times{m}$
\begin{alignat}{5}
y_{1} &=& &x_n- \frac{f(x_n)}{f^\prime(x_n)}, \nonumber \\
y_{2} &=& &x- \frac{f(x_n)}{f^\prime(x_n)}\left[1+\frac{f(y_1)}{f(x_n)}\left(1+2\frac{f(y_1)}{f(x_n)}\right)\right], \nonumber \\
y_{3} &=& &x_n - \frac{f(x_n)}{f^\prime(x_n)}\left[1+\frac{f(y_1)}{f(x_n)}\left(1+2\frac{f(y_1)}{f(x_n)}\right)+\frac{f(y_2)}{f(x_n)}\left(1+2\frac{f(y_1)}{f(x_n)}\right)\right], \nonumber \\
y_{4} &=& &x_n- \frac{f(x_n)}{f^\prime(x_n)}\left[1+\frac{f(y_1)}{f(x_n)}\left(1+2\frac{f(y_1)}{f(x_n)}\right)+\frac{f(y_2)}{f(x_n)}\left(1+2\frac{f(y_1)}{f(x_n)}\right)+\frac{f(y_3)}{f(x)}\left(1+2\frac{f(y_1)}{f(x_n)}\right)\right], \nonumber \\
&\vdots& &\nonumber\\
y_{m-1} &=& &x_n- \frac{f(x_n)}{f^\prime(x_n)}\left[1+\frac{f(y_1)}{f(x_n)}\left(1+2\frac{f(y_1)}{f(x_n)}\right)+\cdots+\frac{f(y_{m-2})}{f(x_n)}\left(1+2\frac{f(y_1)}{f(x_n)}\right)\right], \nonumber \\
x_{n+1} &=& & x_n- \frac{f(x_n)}{f^\prime(x_n)}\left[1+\frac{f(y_1)}{f(x_n)}\left(1+2\frac{f(y_1)}{f(x_n)}\right)+\cdots+\frac{f(y_{m-1})}{f(x_n)}\left(1+2\frac{f(y_1)}{f(x_n)}\right)\right]. \label{eq:scheme}
\end{alignat}
It may be notice that a $2\times{m}$ order method formed by the above scheme will require $m$ functions and one deravitive evaluation per iterative step. A C$^{++}$ implementation of the above scheme \eqref{eq:scheme} is presented in the Listing \ref{magic}.
\lstset{backgroundcolor=\color{white},emph={function,firstderivative},emphstyle={\color{red}\bf}}
\begin{lstlisting}[float,caption=C++ implementation.,label=magic,frame=shadowbox,rulesepcolor=\color{blue},framexleftmargin=5mm,framerule=0pt]
int main(){
    mp::mp_init(2005);
    mp_real tol =  1.0e-300; 
    unsigned int w = 5;
    unsigned int maxitr = 1000;
    std::vector<mp_real> x(maxitr,0.0);
    x[0] = -0.5;
    mp_real err1 = 100.0;
    mp_real err2 = 100.0;
    mp_real gamma = -1.0;
    unsigned int m = 1;
    unsigned i = 0;
    unsigned int conv = 0;
    while(err1 > tol || err2 > tol){
	if(i >= maxitr){
	    break;
	}
	x[i+1] = x[i] - function(x[i])/firstderivative(x[i]);
	mp_real x1 = x[i+1];
	mp_real x0 = x[i];
	for(unsigned int k = 1 ; k < m ; ++k){
	    mp_real y = x[i+1];
	    x[i+1] = x[i+1] - function(x[i+1])
	              /firstderivative(x0)*(1.0 
	                    + 2.0 * function(x1)/function(x0));
	    err1 = abs(y-x[i+1]);
	    err2 = abs(function(x[i+1]));
	    if(err1 < tol & err2 < tol){
		std::cout << "err1 = " << scientific 
			  << setw(2*w) << err1 
			  << "    err2 = " <<setw(2*w) 
		          << err2 <<  std::endl;
		conv = 1;
		break;
	    }
	}
	if(conv) break;
	err1 = abs(x[i+1] - x[i]);
	err2 = abs(function(x[i+1]));
	std::cout << setw(w)  << "itr.=" << i << scientific 
	    	  << "   " << setw(2*w) << x[i+1] 
		  << std::endl;
	++i;
    }
    mp_real rho;
    for (unsigned int n = 1; n < i; ++n){
	rho = log(abs( (x[n+1] - gamma)/(x[n]-gamma))) 
	            /log(abs( (x[n] - gamma)/(x[n-1]-gamma)));
	std::cout << "rho = " << rho << std::endl;
    }
    mp::mp_finalize();
    return 1;
}
\end{lstlisting}
\section{Numerical work}
The order of convergence $\xi$ of an iterative method is defined as follows
\begin{equation}
\lim_{n\to\infty}\dfrac{\vert{e_{n+1}}\vert}{\vert{e_{n}}\vert^\xi} = c \neq{0}.
\end{equation}
Here, $e_n$ is the error after $n$ iterations of the method. Then the approximate value of the computational order of convergence (COC) $\rho$ \citep{chun:numer} can be find using the formula
$$\rho \approx \frac{\textrm{Log}{\vert{({x_{n+1}-\gamma})/({x_n-\gamma})}\vert}}{\textrm{Log}{\vert{({x_{n}-\gamma})/({x_{n-1}-\gamma})}\vert}}.$$

All the computations reported here are done in the programming language C$^{++}$. For numerical precision, we are using ARPREC\citep{arprec}. The ARPREC package supports arbitrarily high level of numeric precision\citep{arprec}. In the program (see the Listing 1) the precision in decimal digits is set at 2005 with the command mp::mp init(2005)\citep{arprec}. We have implemented the presented technique in the C$^{++}$ language. Listing \ref{magic} presents the main part of our implementation.

For convergence, it is required that the distance of two consecutive approximations $(\vert{x_{n+1} - x_{n}}\vert)$ be less than $\epsilon$. And, the absolute value of the function $(\vert{f(x_n)}\vert)$ also referred to as residual be less than $\epsilon$. Apart from the convergence criteria, our algorithm also uses maximum allowed iterations as stopping criterion. Thus our algorithm stops if (i) $\vert{x_{n+1} - x_{n}}\vert< \epsilon$ (ii) $\vert{f(x_n)}\vert< \epsilon$ (iii) $\textrm{itr} > \textrm{maxitr}$. Here, $\epsilon = 1\times{10}^{-300}$, $\textrm{itr}$ is the iteration counter for the algorithm and $\textrm{maxitr}=100$. See the C$^{++}$ algorithm presented in the Listing \ref{magic}. The algorithm is tested for the following functions \citep{eight1}
\begin{alignat}{5}
f_1{(x)} &=& &x^5 + x^4 +4 x^2 -15, \quad &\gamma &\approx& 1.347, \nonumber \\
f_2{(x)} &=& &\sin{x} - x/3,        \quad &\gamma &\approx& 2.278, \nonumber \\
f_3{(x)} &=& &10x\,e^{-x^2}-1, \quad &\gamma &\approx& 1.679, \nonumber \\
f_4{(x)} &=& &\cos{x} - x,        \quad &\gamma &\approx& 0.739, \nonumber \\
f_5{(x)} &=& &e^{-x^2+x+2} - 1,        \quad &\gamma &\approx& -1.000, \nonumber \\
f_6{(x)} &=& & e^{-x}+\cos{x},        \quad &\gamma &\approx& 1.746, \nonumber \\
f_7{(x)} &=& & \textrm{Log}{(x^2+x+2)} -x + 1,        \quad &\gamma &\approx& 4.152, \nonumber \\
f_8{(x)} &=& & \textrm{arcsin}{(x^2-1)}-x/2+1,        \quad &\gamma &\approx& 0.5948. \nonumber 
\end{alignat}
We run the algorithm shown in the Listing \ref{magic} for four values of $m$: $m=1,2,3,4$. Here, $m=1$ corresponds to the classical Newton method. We choose the same initial guess as found in the article \citep{eight1}. So the reader may find it easier to compare performance of various methods. Table \ref{table:1} reports outcome of our numerical work. The table \ref{table:1} reports (iterations required, number of function evaluations needed, COC during second last iteration) for the Newton method ($m=1$), fourth order iterative method ($m=2$), sixth order iterative method ($m=3$) and eight order iterative method ($m=4$). Computational order of convergence reported in the Table \ref{table:1} was observed during the second last iteration.
\begin{table}[!ht]
\centering
\begin{tabular}{|c|c||cccc|}
	\hline
\textbf{\em $f(x)$} & $x_0$ &  NM$(m=1)$    &   $m=2$  & $m=3$ &  $m=4$     \\
	\hline \hline
$f_1(x)$   & $1.6$  & $(9,18,2)$& $(4,12,4)$ & $(2,{\bf{8}},5.66)$ & $(2,10,7.6)$ \\
$f_2(x)$   & $2.0$ & $(23,46,1)$ & $(10,30,1)$ & $(7,{\bf{21}},1)$ & $(6,30,1)$  \\
$f_3(x)$   & $1.8$ & $(10,20,2)$ & $(4,{\bf{12}},3.99)$ & $(3,{\bf{12}},6.21)$ & $(3,15,8.22)$ \\
$f_4(x)$   & $1.0$ & $(9,18,2)$ & $(4,{\bf{12}},3.99)$ & $(3,{\bf{12}},5.90)$ & $(3,15,8.10)$ \\	
$f_5(x)$   & $-0.5$ & $(11,22,2)$ & $(5,{\bf{15}},3.99)$ & $(4,16,5.99)$ & $(3,{\bf{15}},6.75)$ \\ 
$f_6(x)$   & $2.0$ & $(9,18,2)$ & $(4,12,3.99)$ & $(3,12,5.99)$ & $(2,{\bf{10}},8.10)$ \\
$f_7(x)$   &  $3.2$ & $(10,20,2)$ & $(4,{\bf{12}},3.99)$ & $(3,{\bf{12}},6.19)$ & $(3,15,8.19)$ \\
$f_8(x)$   & $1.0$ & $(10,20,2)$ & $(4,{\bf{12}},4.01)$ & $(3,{\bf{12}},6.35)$ & $(3,15,8.36)$ \\\hline
\end{tabular}
\caption{(iterations, number of function evaluations, COC) for the Newton method ($m=1$), fourth order method ($m=2$), sixth order method ($m=3$) and eight order iterative method ($m=4$).}\label{table:1}
\end{table}
Following important observations were made during numerical experimentations
\begin{enumerate}
\item In the Table \ref{table:1}, the methods which require least number of function evaluations for convergence are marked in bold.
\item As reported in the Table \ref{table:1}, for the function $f_2(x)$, COC (during the second last iteration) is same for $m=1,2,3,4$. While COC for the second, third and fourth iterations is reported in the Table \ref{table:2}. 
\item From the Table \ref{table:1}, we notice that for the functions $f_3(x)$, $f_4(x)$, $f_7(x)$ and $f_8(x)$ the sixth order ($m=3$) and eight ($m=4$) methods requirs same number of iterative steps. Table \ref{table:3} reports residual ($\vert{f(x_n)}\vert$) during the last iterative step.
\end{enumerate}

\begin{table}[!ht]
\centering
\begin{tabular}{|c||cccc|}
	\hline
\textbf{$f(x)$}  &  $m=1$    &   $m=2$  & $m=3$ & $m=4$ \\
	\hline \hline
$f_2(x)$   & $10^{-435}$  & $10^{-872}$ & $10^{-671}$ & $10^{-1522}$ \\
	\hline
\end{tabular}
\caption{Computational order of convergence (COC) at the second iteration.}\label{table:2}
\end{table}

\begin{table}[!ht]
\centering
\begin{tabular}{|c||cccc|}
	\hline
\textbf{$f(x)$}  &  $m=1$    &   $m=2$  & $m=3$ & $m=4$ \\
	\hline \hline
$f_7(x)$   & $10^{-435}$  & $10^{-872}$ & $10^{-671}$ & $10^{-1522}$ \\
$f_8(x)$   & $10^{-347}$  & $10^{-744}$ & $10^{-800}$ & $10^{-1302}$ \\
	\hline
\end{tabular}
\caption{Residual ($\vert{f(x_n)}\vert$)}\label{table:3}
\end{table}



\begin{thebibliography}{0}
\bibitem{fourth1}
J. Kou, Y. Li, X. Wang, A composite fourth-order iterative method
for solving non-linear equations, Appl. Math. Comput.  184 (2007)
471-475.

\bibitem{fourth2}
C. Chun, Y. Ham,
\newblock{A one-parameter fourth-order family of iterative methods for nonlinear equations}, Appl. Math. Comput. 189 (2007) 610-614.

\bibitem{fourth3}
J.R. Sharma and R.K. Goyal, Fourth-order derivative-free methods for
solving non-linear equations, Int. J. Comput. Math. 83 (2006)
 101--106.

\bibitem{fourth4}
C. Chun, Some fourth-order iterative methods for solving nonlinear
equations,  Appl. Math. Comput. 195 (2008) 454-459.

\bibitem{fourth5}
M.A. Noor, F. Ahmad, Fourth-order convergent iterative method for
nonlinear equation, Appl. Math. Comput. (2006),
doi:10.1016/j.amc.2006.04.068.

\bibitem{fourth6}
J. Kou, Y. Li, X. Wang, Fourth-order iterative methods free from
second derivative, Appl. Math. Comput. (2006),
doi:10.1016/j.amc.2006.05.189.

\bibitem{fourth7}
R. King, A family of fourth-order methods for nonlinear equations,
SIAM J. Numer. Anal. 10 (1973) 876--879.

\bibitem{fourth8}
C. Chun, A family of composite fourth-order iterative methods for
solving nonlinear equations, Appl. Math. Comput.187 (2007) 951-956.


\bibitem{fourth9}
I.K. Argyros, D. Chen and Q. Qian, The Jarratt method in Banach
space setting, J. Comput. Appl. Math. 51 (1994), 1-3.

\bibitem{fourth10}
C. Chun, Y. Ham, Some fourth-order modifications of Newton's method,
Appl. Math. Comput. 197 (2008) 654-658.

\bibitem{fourth11}
J. Kou, Y. Li, X. Wang, Fourth-order iterative methods free from
second derivative, Appl. Math. Comput. 184 (2007) 880-885.

\bibitem{fourth12}
A. K. Maheshwari, A fourth order iterative method for solving
nonlinear equations, Appl. Math. Comput. In Press, Accepted
Manuscript, Available online 28 January 2009.

\bibitem{xinlog}
X. Feng, Y. He, High order iterative methods without derivatives for
solving nonlinear equations, Appl. Math. Comput. 186 (2007)
1617-1623.

\bibitem{abbasbandy}
S. Abbasbandy, Modified homotopy perturbation method for nonlinear
equations and comparison with Adomian decomposition method, Appl.
Math. Comput. 172 (2006) 431-438.

\bibitem{thirdorder1}
S. Weerakoon and T.G.I. Fernando, A Variant of Newton's Method with
Accelerated Third-Order Convergence, Appl. Math. Lett. 13 (2000)
87-93.

\bibitem{thirdorder2}
H.H. Homeier, On Newton-type methods with cubic convergence, J.
Comput. Appl. Math. 176 (2005) 425-432.

\bibitem{potra1}
F.A. Potra and V. Pt\'{a}k, Nondiscrete induction and iterative
processes, Research Notes in Mathematics vol. 103, Pitman, Boston
(1984).

\bibitem{thirdorder3}
J. Kou, Y. Li, X. Wang, A modification of Newton method with
third-order convergence, Appl. Math. Comput. 181 (2006) 1106-1111.

\bibitem{khattri2}
S. K. Khattri, Altered Jacobian Newton Iterative Method for
Nonlinear Elliptic Problems. IAENG International Journal of Applied
Mathematics, 38:3, IJAM\_38\_3\_01, 2008.

\bibitem{thirdorder4}
M. Frontini and E. Sormani, Some variant of Newton's method with
third-order convergence, Appl. Math. Comput. 140  (2003) 419--426.

\bibitem{thirdorder5}
A.Y. \"{O}zban, Some new variants of Newton's method, Appl. Math.
Lett. 17 (2004) 677--682.

\bibitem{khattri_1}
S. K. Khattri, Newton-Krylov Algorithm with Adaptive Error
Correction for the Poisson-Boltzmann Equation, MATCH Commun. Math.
Comput. Chem.  56 (2006) 197--208.

\bibitem{example}
C. Chun, A geometric construction of iterative functions of order
three to solve nonlinear equations, Comput. Math. Appl.  53  (2007)
972-976.


\bibitem{ostro}
A.M. Ostrowski, Solutions of Equations and System Equations, Academic Press, New Youk, 1960.

\bibitem{traub}
J. F. Traub. Iterative methods for the solution of equations, Chelsea Publishing Company, New Youk, 1977.

\bibitem{jarrat}
I. K. argyros, D. Chen, Q. Qian, The Jarrat method in Banach space setting, J. Computing. Appl. Math. 51 (1994).

\bibitem{chun:numer}
Changbum Chun. Construction of Newton-like iteration methods for solving nonlinear equations. Numerische Mathematik. Volume 104, Number 3 / September, 2006.



\bibitem{six1}
C. Chun, Y. Ham. Some sixth-order variants of Ostrowski root-finding methods. Appl. Math. Comput. 193, 2003

\bibitem{six2} 
H. Ren, Q. Wu and W. Bi. New variants of Jarratts method with sixth-order convergence. 	Numerical Algorithms, Volume 52, Number 4, December, 2009.


\bibitem{six3}
J.R. Sharma, and R.K. Guha.
A family of modified Ostrowski methods with accelerated sixth order convergence. Appl. Math. Comput. 190, 2007. 


\bibitem{eight1}
Weihong Bi, Hongmin Ren and Qingbiao Wu. Three-step iterative methods with eighth-order convergence for solving nonlinear equations. Journal of Computational and Applied Mathematics, Volume 225, Issue 1, 1 March 2009, Pages 105-112.

\bibitem{arprec}
ARPREC. C++/Fortran-90 arbitrary precision package. Available at \url{http://crd.lbl.gov/~dhbailey/mpdist/}.

\end{thebibliography}
\end{document}